\documentclass[12pt]{amsart}
\usepackage{amsmath,amssymb,amsthm}

\newcommand{\R}{\mathbb R}
\newcommand{\Sp}{\mathbb S}

\renewcommand{\span}{\mathrm{span}}

\newtheorem{thm}{Theorem}[section]

\theoremstyle{definition}

\theoremstyle{remark}
\newtheorem{rem}[thm]{Remark}
\newcommand{\ds}{\displaystyle}

\begin{document}

\title[MINIMAL HYPERSURFACES OF TYPE NUMBER TWO]
{MINIMAL SURFACES IN THE THREE-DIMENSIONAL SPHERE AND MINIMAL HYPERSURFACES OF TYPE NUMBER TWO}%

\author{GEORGI GANCHEV}%
\address{Bulgarian Academy of Sciences, Institute of Mathematics and
Informatics, Acad. G. Bonchev Str. bl. 8, 1113 Sofia, Bulgaria}%
\email{ganchev@math.bas.bg}%

\subjclass[2000]{Primary 53A07, Secondary 53A10}%

\keywords{Minimal surfaces in the three-dimensional sphere; canonical parameters;
bi-umbilical hypersurfaces of type number two; minimal hypersurfaces of type number two}%

\begin{abstract}
We introduce canonical principal parameters on any strongly regular minimal surface in the
three dimensional sphere and prove that any such a surface is determined up to a motion
by its normal curvature function satisfying the Sinh-Poisson equation.
We obtain a classification theorem for bi-umbilical hypersurfaces of type number two.
We prove that any minimal hypersurface of type number two with involutive distribution
is generated by a minimal surface in the three-dimensional Euclidean space, or in the
three dimensional sphere. Thus we prove that the theory of minimal hypersurfaces of
type number two with involutive distribution is locally equivalent to the theory of
minimal surfaces in the three dimensional Euclidean space or in the three-dimensional
sphere.
\end{abstract}

\maketitle
\section{Introduction}
The aim of this paper is to show the deep relation between the theory of minimal
hypersurfaces of type number two in arbitrary dimension and the theory of minimal
surfaces in the three-dimensional Euclidean space ${\R}^3$ or in the three-dimensional
Euclidean sphere ${\Sp}^3$.

Section 1 is devoted to the invariant theory of minimal surfaces in ${\Sp}^3$.

Let $\mathcal M: \, z=l(u,v), \; (u,v) \in {\mathcal D}$ be a regular
surface in ${\Sp}^3$ parameterized by principal parameters $(u,v)$ and consider
the following four invariant functions: the principle normal curvatures $\nu_1, \, \nu_2$;
the principal geodesic curvatures (the geodesic curvatures of the principal lines)
$\gamma_1, \, \gamma_2$. We introduce the class of strongly regular surfaces by
the condition
$$(\nu_1-\nu_2)\, \gamma_1\, \gamma_2 \neq 0.$$

In Subsection 2.2 we prove the Bonnet type Theorem \ref{T:2.1} for strongly regular
surfaces in terms of the four invariants $\nu_1, \, \nu_2, \, \gamma_1, \, \gamma_2$.

The main obstacle to formulate the Bonnet type fundamental theorem for surfaces
similarly to the fundamental theorem for the curves is the lack of natural parameters
in the theory of surfaces. In \cite{GM} we have shown that the class of
Weingarten surfaces in ${\R}^3$ admits geometrically determined parameters,
which we called \emph{canonical} parameters. Here we show that any minimal
strongly regular surface in ${\Sp}^3$ admits geometrically determined canonical
parameters.

Our main result for minimal strongly regular surfaces in ${\Sp}^3$ is Theorem \ref{T:2.2}:
\vskip 1mm
\emph{Any solution $\nu>0$ to the Sinh-Poisson equation
$$\Delta \ln \nu=2\frac{1-\nu^2}{\nu} $$
satisfying the condition $\nu_u \, \nu_v \neq 0$ determines uniquely $($up to a motion
in ${\Sp}^3$$)$ a minimal strongly regular surface with invariants
$$\nu_1=\nu, \quad \nu_2=-\nu, \quad \gamma_1=(\sqrt {\nu})_v,  \quad
\gamma_2= -(\sqrt{\nu})_u.$$
Furthermore $(u,v)$ are canonical parameters.}
\vskip 1mm
This result allows us to introduce a family $\{\mathcal M_t\}$ of minimal surfaces
associated with a given minimal strongly regular surface $\mathcal M$. The minimal surfaces
of the family $\{\mathcal M_t\}$ are isometric to the surface $\mathcal M$.

Section 2 is devoted to minimal hypersurfaces of type number two.

Let $({\mathcal M}^n,g), \; n\geq 3$ be a regular hypersurface in the $(n+1)$-dimensional
Euclidean space ${\R}^{n+1}$. If the shape operator of the hypersurface ${\mathcal  M}^n$
has two eigenvalues $\nu_1, \, \nu_2$ different from zero and the other $(n-2)$
eigenvalues are zero at each point, then ${\mathcal M}^n$ is said to be of {\it type
number two}. Thus, the hypersurfaces of type number two are characterized in
terms of their fundamental form $h$ as follows:
$$h=\nu_1\, \eta_1 \otimes \eta_1+ \nu_2 \, \eta_2 \otimes \eta_2, \quad
\nu_1\,\nu_2 \neq 0,$$
where $\eta_1, \, \eta_2$ are unit one-forms.

These hypersurfaces considered as Riemannian manifolds are semi-symmetric spaces
foliated by Euclidean leaves of codimension two, i.e. they are \emph{foliated
semi-symmetric spaces} \cite{S}. Conversely, any foliated semi-symmetric hypersurface
in ${\R}^{n+1}$ is a hypersurface of type number two. Foliated semi-symmetric spaces
have been studied in \cite{BKV} with respect to their metrics under the name Riemannian
manifolds of {\it conullity two}. Thus, the hypersurfaces of type number two
can also be considered as hypersurfaces of conullity two.

In \cite{GM2} we considered the hypersurfaces of type number two as the envelope of a
two parameter system of hyperplanes. Then, a hypersurface ${\mathcal M}^n$ of type number
two can be considered as a two-parameter system of planes of codimension three
$\{E^{n-2}(u,v)\}, \; (u,v) \in \mathcal D$ with the property that
the tangent hyperplane to ${\mathcal M}^n$ at any point of an arbitrary generator
$E^{n-2}$ is one and and the same. Briefly, the hypersurfaces of type number two are
exactly the "developable" two parameter systems of planes of codimension three.

The unit eigenvector fields $X$ and $Y$ corresponding to $\nu_1$ and $\nu_2$ determine a
two-dimensional distribution $\Delta=\span\{X,Y\}$, which plays an essential role in the
geometry of hypersurfaces of type number two. We denote by $\mathcal K_0$ the class of
hypersurfaces of type number two, whose distribution $\Delta$ is involutive.

In \cite{GM1}  we proved that a hypersurface ${\mathcal M}^n$ of type number two is
in the class $\mathcal K_0$, then any two-dimensional integral surface $\mathcal M^2$ of its
distribution $\Delta$ has flat normal connection. Conversely, for any surface $\mathcal M^2$
with flat normal connection, using the set of its parallel surfaces, we gave a natural
geometric construction of a family of hypersurfaces of type number two belonging to
the class $\mathcal K_0$. Thus the local differential geometry of hypersurfaces in
the class $\mathcal K_0$ is equivalent to the local differential geometry of the
surfaces with flat normal connection. Generally speaking, the theory of the surfaces
with flat normal connection can be considered as a model of the theory of the
hypersurfaces in the class $\mathcal K_0$. The aim of this paper is to realize
this correspondence for the the minimal hypersurfaces from the class $\mathcal K_0$.

We note that the theory of hypersurfaces of type number two carries some features of the
usual theory of surfaces in the Euclidean space ${\R}^3$.

A point of a hypersurface ${\mathcal M}^n$ of type number two is said to be bi-umbilical
\cite{GM2} if $\nu_1=\nu_2\neq 0$ at this point. This notion corresponds to the notion
of an umbilical point in Euclidean differential geometry. In Section 2 we obtain the
classification Theorem \ref{T:3.1} for bi-umbilical hypersurfaces of type number two
(corresponding to the classical result for umbilical surfaces). Our scheme is the following:
\vskip 1mm
We prove that the integral surfaces of the distribution $\Delta$ of any bi-umbilical
hypersurface of type number two lie on two-dimensional spheres. Conversely, any two-dimensional
sphere generates a family of bi-umbilical hypersurfaces of type number two.
\vskip 1mm
Any bi-umbilical hypersurface of type number two is determined by a unit vector function
$l(u,v) \in {\R}^{n+1}$ and a scalar function $r(u,v)$ of two variables satisfying the system
of partial differential equations \cite{GM2}:
$$\begin{array}{cl}
l_{uu} - l_{vv} & = \displaystyle{\frac{E_u}{E}\,l_u - \frac{E_v}{E}\,l_v},\\
[2mm]
2 l_{uv} & = \displaystyle{\frac{E_v}{E}\,l_u + \frac{E_u}{E}\,l_v},\\
[3mm]
r_{uu} - r_{vv}  & = \displaystyle{\frac{E_u}{E}\,r_u - \frac{E_v}{E}\,r_v},\\
[2mm]
2 r_{uv} & = \displaystyle{\frac{E_v}{E}\,r_u + \frac{E_u}{E}\,r_v},
\end{array}
\qquad l_u^2=l_v^2=E, \quad l_u \,l_v=0.$$
\vskip 1mm
Theorem \ref{T:3.1} implies that the solutions of the above system can be found explicitly.
\vskip 1mm
The classification of minimal ruled hypersurfaces (helicoids) in ${\R}^{n+1}$ was treated
in a series of papers (e.g. \cite{FG}, \cite{A2}) and was completed in 1981 \cite{A1}. It occurred
that essential helicoids exist only in ${\R}^3$ (the usual helicoids) and in ${\R}^4$
(second type helicoids). While the usual helicoid ${\mathcal M}^2$ in ${\R}^3$ is a complete surface,
the three-dimensional second type helicoid ${\mathcal M}^3$ in ${\R}^4$ has one singular point. All
minimal ruled hypersurfaces in an arbitrary dimension are generated by these helicoids.

We note that the ruled hypersurfaces are in the class $\mathcal K_0$. In the present paper we
obtain a geometric description of the minimal hypersurfaces in the class $\mathcal K_0$.
This is obtained in Theorem \ref{T:4.1}. Our scheme is the following:
\vskip 1mm
We prove that any integral surface of the distribution $\Delta$ of a minimal hypersurface
of the class $\mathcal K_0$ is a minimal surface in the three-dimensional Euclidean space
${\R}^3$, or a minimal surface in the three-dimensional Euclidean sphere ${\Sp}^3$.
Conversely, any minimal surface in $\R^3$ or in ${\Sp}^3$ generates a minimal
hypersurface of the class $\mathcal K_0$.
\vskip 2mm
Briefly, Theorem \ref{T:4.1} means that the model theory of minimal hypersurfaces in the class
$\mathcal K_0$ is the theory of minimal surfaces in ${\R}^3$ and in ${\Sp}^3$.

Any minimal hypersurface of the class $\mathcal K_0$ is determined by a unit vector function
$l(u,v) \in {\R}^{n+1}$ and a scalar function $r(u,v)$ of two variables satisfying the system
of partial differential equations \cite{GM2}:
$$\begin{array}{l}
l_{uu} + l_{vv} + 2E\,l = 0,\\
[2mm]
r_{uu} + r_{vv} + 2E\,r = 0,
\end{array} \qquad l_u^2=l_v^2=E, \quad l_u \, l_v=0.$$
Theorem \ref{T:4.1} implies that any example of a minimal surface in ${\R}^3$ or in ${\Sp}^3$
generates a family of solutions to the above system.
\section{An invariant theory of minimal surfaces in the three-dimensional sphere}
\subsection{Strongly regular surfaces in ${\Sp}^3$}
We consider the three dimensional sphere ${\Sp}^3$ as a unit hyper-sphere centered at the origin
in the Euclidean space ${\R}^4$.  The unit normal vector field to ${\Sp}^3$ is denoted by
$l$ and the flat Levi-Civita connection of the standard metric in ${\R}^4$ is denoted by
$\nabla$.

Let $\mathcal M: \, z=l(u,v), \; (u,v) \in {\mathcal D}$ be a surface
in ${\Sp}^3$. We denote by $N$ the unit normal vector field to ${\mathcal M}$ and by
$E, F, G; \; e, f, g$ - the coefficients of the first and the second
fundamental forms, respectively.

We suppose that the surface has no umbilical points and the principal lines on
$\mathcal M$ form a parametric net, i. e.
$$F(u,v)=f(u,v)=0, \quad (u,v) \in \mathcal D.$$
Then the principal curvatures $\nu_1, \nu_2$ and the principal geodesic curvatures
(geodesic curvatures of the principal lines) $\gamma_1, \gamma_2$ are given by
$$\nu_1=\frac{e}{E}, \quad \nu_2=\frac{g}{G}; \qquad
\gamma_1=-\frac{E_v}{2E\sqrt G}, \quad \gamma_2= \frac{G_u}{2G\sqrt E}. \leqno(2.1)$$

We consider the canonical tangential frame field $\{X, Y\}$ determined by
$$X:=\frac{z_u}{\sqrt E}, \qquad Y:=\frac{z_v}{\sqrt G}.$$

The following Frenet type formulas for the orthonormal frame field $XYNl$
associated with every point of ${\mathcal M}$ are valid

$$
\begin{array}{clllll}
\nabla_X \,X &= & & \; \; \; \gamma_1 \,Y & +\nu_1 \, N & + l,\\
[2mm]
\nabla_X \,Y &= &-\gamma_1 \, X, & \\
[2mm]
\nabla_X \,N &= &-\nu_1 \, X, & \\
[2mm]
\nabla_X \,l &= &\;\;\;-\,X & \\
[5mm]
\nabla_Y \,X &= & & \; \; \; \gamma_2 \,Y, \\
[2mm]
\nabla_Y \,Y &= &-\gamma_2 \,X & &+ \nu_2 \,N & + l,\\
[2mm]
\nabla_Y \,N &= & & - \nu_2 \,Y,\\
[2mm]
\nabla_Y \,l &=&& \;\;\;-\,Y.
\end{array}
\leqno(2.2)$$
\vskip 2mm

The Codazzi equations have the following form
$$\gamma_1=\frac{Y(\nu_1)}{\nu_1-\nu_2}= \frac{(\nu_1)_v}{\sqrt G(\nu_1-\nu_2)},
\qquad
\gamma_2=\frac{X(\nu_2)}{\nu_1-\nu_2}=
\frac{(\nu_2)_u}{\sqrt E(\nu_1-\nu_2)} \leqno(2.3)$$
and the Gauss equation can be written as follows:
$$Y(\gamma_1)-X(\gamma_2)-(\gamma_1^2+\gamma_2^2)=1+\nu_1\nu_2,$$
or
$$\frac{(\gamma_1)_v}{\sqrt G}-\frac{(\gamma_2)_u}{\sqrt E}-
(\gamma_1^2+\gamma_2^2)=1+\nu_1\nu_2.\leqno(2.4)$$

\begin{rem}
The mean curvature of $\mathcal M$ in ${\Sp}^3$ is the function $(\nu_1+\nu_2)/2$,
while the Gauss curvature (the sectional curvature) of ${\mathcal M}$ is $K=1+\nu_1\,\nu_2$.
\end{rem}

A surface ${\mathcal M}: \; z=l(u,v), \; (u,v)\in \mathcal D$
parameterized with principal parameters is said to be \emph{strongly regular} if
(cf \cite{GM})
$$\gamma_1(u,v)\gamma_2(u,v) \neq 0, \quad (u,v)\in\mathcal D.$$

Since $$\gamma_1 \gamma_2 \neq 0 \; \iff \; (\nu_1)_v (\nu_2)_u \neq 0,$$
then the following formulas
$$\sqrt E=\frac{(\nu_2)_u}{\gamma_2(\nu_1-\nu_2)} >0, \quad
\sqrt G=\frac{(\nu_1)_v}{\gamma_1 (\nu_1-\nu_2)}>0\, .\leqno(2.5)$$
are valid for any strongly regular surface. Because of (2.5) formulas (2.2) become

$$\begin{array}{lcccc}
X_u  = &  & \displaystyle{\frac{\gamma_1 \, (\nu_2)_u} {\gamma_2(\nu_1-\nu_2)}\, Y}
&+ \, \displaystyle{\frac{\nu_1 \, (\nu_2)_u}{\gamma_2(\nu_1-\nu_2)}\, N}&
+\displaystyle{\frac{(\nu_2)_u}{\gamma_2(\nu_1-\nu_2)}\, l},\\
[4mm]
Y_u  = & - \, \displaystyle{\frac{\gamma_1 \, (\nu_2)_u}
{\gamma_2(\nu_1-\nu_2)}\, X}, & & &\\
[4mm]
N_u  = & - \, \displaystyle{\frac{\nu_1 \, (\nu_2)_u}
{\gamma_2(\nu_1-\nu_2)}}\, X,  & & &\\
[3mm]
l_u = & - \, \displaystyle{\frac{(\nu_2)_u}{\gamma_2(\nu_1-\nu_2)}\, X;}\\
[5mm]
X_v = & &\displaystyle{\frac{\gamma_2 \, (\nu_1)_v}
{\gamma_1(\nu_1-\nu_2)}\, Y,} &\\
[4mm]
Y_v  = &- \, \displaystyle{\frac{\gamma_2 \, (\nu_1)_v}
{\gamma_1(\nu_1-\nu_2)}}\, X & &
+ \, \displaystyle{\frac{\nu_2 \,(\nu_1)_v}{\gamma_1
(\nu_1-\nu_2)}}\, N &+\displaystyle{\frac{(\nu_1)_v}{\gamma_1(\nu_1-\nu_2)}\, l},\\
[4mm]
N_v  = &  & - \, \displaystyle{\frac{\nu_2 \, (\nu_1)_v}
{\gamma_1(\nu_1-\nu_2)}}\,Y, &\\
[4mm]
l_v = & & - \, \displaystyle{\frac{(\nu_1)_v}{\gamma_1(\nu_1-\nu_2)}\, Y.}
\end{array}
\leqno(2.6)$$
Then the fundamental theorem for strongly regular surfaces states in terms of the
four invariants $\nu_1, \, \nu_2, \, \gamma_1, \, \gamma_2$ as follows:
\begin{thm}\label{T:2.1}
\textbf{$($Bonnet type fundamental theorem$)$}
Given four functions $\nu_1(u,v)$, \, $\nu_2(u,v)$, \, $\gamma_1(u,v)$,
\, $\gamma_2(u,v)$, \; $(u,v)\in \mathcal D$ satisfying the following
conditions:
$$\begin{array}{ll}
1) & \nu_1-\nu_2>0, \quad \gamma_1 \, (\nu_1)_v >0, \quad
\gamma_2 \, (\nu_2)_u > 0; \\
[4mm]
2.1) & \displaystyle{\left(\ln\frac{(\nu_1)_v}{\gamma_1}\right)_u=
\frac{(\nu_1)_u}{\nu_1-\nu_2},}
\qquad
\displaystyle{\left(\ln\frac{(\nu_2)_u}{\gamma_2}\right)_v=
-\frac{(\nu_2)_v}{\nu_1-\nu_2};}\\
[5mm]
2.2) & \displaystyle{\frac{\nu_1-\nu_2}{2}\left(\frac{(\gamma_1^2)_v}
{(\nu_1)_v}-\frac{(\gamma_2^2)_u}{(\nu_2)_u}\right)-
(\gamma_1^2+\gamma_2^2)=1+\nu_1\nu_2,}
\end{array}$$
and an initial right oriented orthonormal frame $l_0X_0Y_0N_0$.

Then there exists a unique strongly regular surface
${\mathcal M}: \; z=l(u, v), \; (u, v) \in \mathcal D_0 \;
((u_0, v_0) \in \mathcal D_0 \subset \mathcal D)$ in ${\Sp}^3$, such that

$(i)$ \; \; $(u, v)$ are principal parameters;

$(ii)$ \, \, $l(u_0, v_0)=l_0, \; X(u_0, v_0)=X_0, \;
Y(u_0, v_0)=Y_0, \; N(u_0,v_0)=N_0$;

$(iii)$ \, the invariants of ${\mathcal M}$ are the given functions
$\nu_1, \, \nu_2, \, \gamma_1, \, \gamma_2.$
\end{thm}
\emph{Proof:} Let $X(u,v), \, Y(u,v), \, N(u,v), \, l(u,v); \; (u,v) \in \mathcal D$ be
four unknown vector valued functions in ${\R}^4$, satisfying the system (2.6).
We write the system (2.6) in the form
$$\left(
\begin{array}{l}
X_u\\
[2mm]
Y_u\\
[2mm]
N_u\\
[2mm]
l_u \end{array} \right)=
A \left(
\begin{array}{l}
X\\
[2mm]
Y\\
[2mm]
N\\
[2mm]
l  \end{array} \right), \qquad
\left(
\begin{array}{l}
X_v\\
[2mm]
Y_v\\
[2mm]
N_v\\
[2mm]
l_v \end{array} \right)=
B \left(
\begin{array}{l}
X\\
[2mm]
Y\\
[2mm]
N\\
[2mm]
l  \end{array} \right),\leqno(2.7)
$$
where $A$ and $B$ are the skew symmetric $4\times4$ matrices determined by (2.6).

The integrability condition of the system (2.7) is given by the equality
$$B_u-A_v=[A,B]. \leqno(2.8)$$

Taking into account (2.6), we find that (2.8) is equivalent to the
conditions 2.1) and 2.2) of the theorem.

Now, let $l_0X_0Y_0N_0$ be an initial orthonormal right oriented
coordinate system at the point $l_0 \in {\Sp}^3$. Applying the theorem
of existence and uniqueness of a solution to (2.6) with initial conditions
$$l(u_0,v_0)=l_0, \quad X(u_0,v_0)=X_0, \quad Y(u_0,v_0)=Y_0, \quad N(u_0,v_0)=N_0,$$
we obtain a uniquely determined solution $l(u,v), \, X(u,v), \, Y(u,v), \, N(u,v); \;
(u,v) \in \mathcal D'$, $(u_0,v_0)\in \mathcal D'\subset \mathcal D$.
\vskip 2mm
Further, we have to prove that $lXYN$ form an orthonormal right oriented frame
field in $\mathcal D'$.

Let $l(l^1,l^2,l^3,l^4), \, X(X^1,X^2,X^3,X^4), \, Y(Y^1,Y^2,Y^3,Y^4), \, N(N^1,N^2,N^3,N^4),$
and set
$$f^{ij}(u,v):=l^il^j+X^iX^j+Y^iY^j+N^iN^j, \quad i,j=1,2,3,4.\leqno(2.9)$$

Differentiating (2.9) with respect to $u$ and $v$ and taking into account (2.6), we obtain
that $f^{ij}(u,v)={\rm const}=\delta_{ij}$, $\delta_{ij}$ being the
Kronecker's deltas. This proves that $lXYN$ form an orthonormal
right oriented frame field at any point $(u,v) \in \mathcal D'$.

Now let us consider the surface $\mathcal M: \; z= l(u,v)$, \, $(u,v)\in \mathcal D'$.

Taking into account that the vector valued functions $X, \, Y, \, N, \, l$ satisfy (2.6),
we find that the functions $E=l_u^2, \; G=l_v^2$ satisfy (2.5), which implies that
the invariants of $\mathcal M$ are the given functions
$\nu_1, \, \nu_2, \, \gamma_1, \, \gamma_2$. \hfill{\qed}
\subsection{ Minimal strongly regular surfaces in ${\Sp}^3$}
Let ${\mathcal M}: \; z=l(u,v), \; (u,v)\in \mathcal D$ be a minimal surface in ${\Sp}^3$,
i.e. $\nu_1+\nu_2=0$. We set $\nu=\nu_1$ and further assume that $\nu>0$. Then we have
$$\nu_1=\nu, \qquad \nu_2=-\nu.$$
Further we assume that ${\mathcal M}$ is a strongly regular surface, i.e.
$$\nu_u\,\nu_v\neq 0.$$
Taking into account the equalities (2.1) and (2.3) we find
$$\gamma_1=-\frac{(\ln \sqrt E)_v}
{\sqrt G}=\frac{(\ln \sqrt \nu)_v}{\sqrt G}, \quad
\gamma_2=\frac{(\ln \sqrt G)_u}{\sqrt E}=-\frac{(\ln \sqrt \nu)_u}{\sqrt E},$$
which imply that $\sqrt{\nu E}$ does not depend on $v$, while $\sqrt{\nu G}$
does not depend on $u$.

Let $(u_0, v_0)$ be a fixed point in $\mathcal D$. We introduce new parameters
in a neighborhood of $(u_0, v_0)$ by the formulas
$$\bar u= \int_{u_0}^u\,\sqrt{\nu\,E}\,du, \quad
\bar v= \int_{v_0}^v\,\sqrt{\nu\,G}\,dv.$$
and call them \emph{canonical parameters} (cf \cite{GM}). It is easy to check that
$$\bar E=\bar G=\frac{1}{\nu}$$
with respect to the canonical parameters $(\bar u, \bar v)$.

Further we assume, that the minimal strongly regular surface
$\mathcal M: \, z=l(u,v), \; (u,v) \in {\mathcal D}$
is parameterized with canonical principal parameters. Then we have
$$E=\frac{1}{\nu}, \quad G=\frac{1}{\nu}, \quad
\gamma_1=(\sqrt{\nu})_v, \quad \gamma_2=-(\sqrt{\nu})_u.\leqno(2.10)$$

Theorem \ref{T:2.1}, applied to minimal strongly regular surfaces parameterized
with canonical principal parameters states as follows:
\begin{thm}\label{T:2.2}
Given a function $\nu (u,v) > 0$
in a neighborhood $\mathcal D$ of $(u_0, v_0)$ with $\nu_u \nu_v \neq 0$,
satisfying the partial differential equation
$$\Delta \ln \nu=2\frac{1-\nu^2}{\nu}
\quad (\Delta {\rm \, - \, Laplace \, operator}) \leqno(2.11)$$
and an initial right oriented orthonormal frame $l_0X_0Y_0N_0$.

Then there exists a unique minimal strongly regular surface
${\mathcal M}: \; z=l(u, v), \; (u, v) \in \mathcal D_0 \;
((u_0, v_0) \in \mathcal D_0 \subset \mathcal D)$ in ${\Sp}^3$, such that

$(i)$ \; \; $(u, v)$ are canonical principal parameters;

$(ii)$ \, \, $l(u_0, v_0)=l_0, \; X(u_0, v_0)=X_0, \;
Y(u_0, v_0)=Y_0, \; N(u_0,v_0)=N_0$;

$(iii)$ \, the invariants of ${\mathcal M}$ are the following functions
$$\nu_1=\nu, \quad \nu_2= -\nu, \quad \gamma_1=(\sqrt {\nu})_v,  \quad
\gamma_2= -(\sqrt{\nu})_u.$$
\end{thm}
\begin{rem}
Under the equalities  (2.10) the integrability  conditions 2.1) and 2.2) in
Theorem \ref{T:2.1} reduce to (2.11). Putting $f=\ln \nu$, the partial differential
equation (2.11) gets the form
$$\Delta f +4\,\sinh f=0 \quad ({\rm sinh-Poisson \, \, equation}).$$
\end{rem}

Thus the sinh-Poisson equation is the natural partial differential equation of
minimal (strongly regular) surfaces in ${\Sp}^3$.

Theorem \ref{T:2.2} gives locally a one-to-one correspondence between minimal strongly
regular surfaces (considered up to a motion in ${\Sp}^3$) and the solutions of the natural
partial differential equation (2.11), satisfying the conditions
$$\nu > 0, \quad \nu_u\nu_y\neq 0. \leqno(2.12)$$
\vskip 2mm
\emph{Associated minimal surfaces to a given minimal strongly regular surface}

A beautiful fact in the theory of minimal surfaces in the three-dimensional Euclidean space
is that any minimal surface generates a family of associated minimal surfaces which are
isometric to the given one.

Using Theorem \ref{T:2.2} we introduce a family of associated minimal surfaces in
the following way.

Let $\mathcal M: \, z=z(u,v), \; (u,v) \in {\mathcal D}$ be a minimal
strongly regular surface in ${\Sp}^3$ parameterized by canonical parameters. Assume
that the domain $\mathcal D$ is a disc centered at the origin of the parametric plane.
If $\nu(u,v)$ is the normal curvature function of $\mathcal M$, then it is a solution
to (2.11) satisfying the conditions (2.12). Putting
$$\nu_t(u,v):=\nu(\cos t \, u-\sin t \, v, \; \sin t \, u+\cos t \, v), \quad t\in[0, 2 \pi),
\quad (u,v)\in \mathcal D,$$
we obtain a family of solutions $\{\nu_t,\;t\in [0,2 \pi)\}$ to (2.11) satisfying
(2.12). Applying Theorem \ref{T:2.2} to any solution $\nu_t$, we obtain a family  of minimal
strongly regular surfaces $\{{\mathcal M}_t\}$. It follows immediately that
$$\begin{array}{l}
E_t(u,v)=E(\cos t \, u-\sin t \, v, \; \sin t \, u+\cos t \, v)\\
[2mm]
G_t(u,v)=G(\cos t \, u-\sin t \, v, \; \sin t \, u+\cos t \, v).
\end{array}$$
i.e. any surface $\mathcal M_t$ is isometric to the the given surface $\mathcal M$. It is
natural to call $\{\mathcal M_t, \,t\in [0,2 \pi)\}$ a family of minimal surfaces
\emph{associated} with a given strongly regular surface $\mathcal M$.
\vskip 1mm
\begin{rem}
The above approach for introducing a family of associated isometric surfaces
can be applied to many other cases. For example, applying the results in \cite{GM}, \cite{GM3},
it gives a family of:

- surfaces of constant mean curvature isometric with a given CMC-surface in ${\R}^3$;

- surfaces of constant Gauss curvature ($K=-1$) isometric with a given surface of  constant Gauss
  curvature ($K=-1$) in ${\R}^3$;

- minimal surfaces isometric with a given minimal surface in ${\R}^4.$
\end{rem}
\section{Minimal hypersurfaces of conullity two}
\subsection{Hypersurfaces of conullity two}
Let $({\mathcal M}^n,g)$ be a regular hypersurface in the Euclidean space $({\R}^{n+1},g)$ with the
induced metric $g$ and shape operator $A$. The standard flat Levi-Civita connection of
$({\R}^{n+1},g)$ is denoted by $\nabla$.

Let the hypersurface ${\mathcal M}^n$ be of \emph{type number two}, i.e. its shape operator $A$ has
two different from zero eigenvalues $\nu_1$ and $\nu_2$ and the other $n-2$ eigenvalues
are zero at each point. We recall that the case $\nu_1=\nu_2=0$ characterizes locally
a hyperplane and the case $\nu_1\neq 0,\; \nu_2=0 $ (hypersurfaces with type number one)
characterizes locally a developable hypersurface. Here we study the case
$$\nu_1(p)\,\nu_2(p)\neq 0, \quad p \in {\mathcal M}^n.$$

Since our considerations are local, we can choose two unit tangent vector fields
$X$ and $Y$, such that
$$AX=-\nu_1\,X, \quad AY=-\nu_2\,Y.$$
Then for any tangent vector $x_0 \perp X,Y$ we have
$$Ax_0=0.$$
We denote by $\Delta$ the distribution  $\span\{X,Y\}$, and by $\Delta^{\perp}$ the
distribution, orthogonal to $\Delta$.

Let us introduce the 1-form $\sigma$ on the distribution $\Delta^{\perp}$ by
the formula
$$\sigma(x_0)=g(\nabla_{x_0}X,Y), \quad x_0\in \Delta^{\perp}.$$

Then the standard Codazzi equations for the shape operator $A$ imply the following formulas
$$\begin{array}{l}
g(\nabla_XX,x_0)=d\,\ln \,\nu_1(x_0),\\
[3mm]
g(\nabla_YY,x_0)=d\,\ln \,\nu_2(x_0),\\
[3mm]
g(\nabla_XY,x_0)=\ds{\frac{\nu_1-\nu_2}{\nu_2}\,\sigma(x_0)},\\
[3mm]
g(\nabla_YX,x_0)=\ds{\frac{\nu_1-\nu_2}{\nu_1}\,\sigma(x_0)},\\
[3mm]
\nabla_{x_0}X=\sigma(x_0)\,Y,\\
[3mm]
\nabla_{x_0}Y=-\sigma(x_0)\,X,\\
[3mm]
(\nu_1-\nu_2)^2\,\gamma_1=(\nu_1-\nu_2)Y(\nu_1),\\
[3mm]
(\nu_1-\nu_2)^2\,\gamma_2=(\nu_1-\nu_2)X(\nu_2),
\end{array}\quad \quad \quad  x_0\in \Delta^{\perp}.\leqno(3.1)$$

Let $N$ be the unit normal vector field to ${\mathcal M}^n$ and $\{e_1,...,e_{n-2}\}$ be an orthonormal
frame field such that $\Delta^{\perp}=\span\{e_1,...,e_{n-2}\}$. Taking into account (3.1),
we write the following Frenet type formulas for the derivatives with respect to $X$ and $Y$:
$$
\begin{array}{lccl}
\nabla_XX= & \gamma_1\,Y & + \, \ds{\sum_i \lambda_i \, e_i} &+\,\nu_1\,N,\\
[4mm]
\nabla_XY=-\gamma_1\,X &  &
+\,\ds{\frac{\nu_1-\nu_2}{\nu_2} \,\sum_i \sigma_i \, e_i,} &\\
[4mm]
\nabla_Xe_i=-\lambda_i\,X &-\ds{\frac{\nu_1-\nu_2}{\nu_2}\,\sigma_i\,Y}, & &\\
[4mm] \nabla_XN=-\nu_1\,X,\\
[8mm]
\nabla_YX= &\gamma_2\,Y &+ \,\ds{\frac{\nu_1-\nu_2}{\nu_1} \,\sum_i \sigma_i \, e_i,} &\\
[4mm]
\nabla_YY=-\gamma_2\,X & & +\,\ds{\sum_i \mu_i \, e_i} &+\,\nu_2\,N,\\
[4mm]
\nabla_Ye_i=-\,\ds{\frac{\nu_1-\nu_2}{\nu_1} \,\sigma_i\,X}& -\mu_i\,Y, & & \\
[4mm]
\nabla_YN=& -\nu_2\,Y, & &\\
\end{array}\leqno(3.2)
$$
where
$$ \gamma_1=g(\nabla_XX,Y), \; \gamma_2=-g(\nabla_YY,X); \; \lambda_i=
\frac{e_i(\nu_1)}{\nu_1}, \; \mu_i=\frac{e_i(\nu_2)}{\nu_2}, \; \sigma_i=
\sigma(e_i); \; i=1,...,n-2.$$

First we shall consider the special class of bi-umbilical hypersurfaces of type number two,
which is the analogue to the umbilical surfaces in ${\R}^3$. In this case the tangent space
at each point of ${\mathcal M}^n$ consists of two orthogonal umbilical distributions
$\Delta$ and $\Delta^{\perp}$.
\subsection{Bi-umbilical hypersurfaces of type number two}
Let us consider the class of hypersurfaces of type number two satisfying the
equality
$$\nu_1(p)=\nu_2(p)=\nu(p), \quad p\in {\mathcal M}^n.$$
Since the shape operator $A$ has equal eigen values at each of the mutually orthogonal
distributions $\Delta$ and $\Delta^{\perp}$, we called these hypersurfaces
\emph{bi-umbilical} \cite{GM2}.

The first property of these hypersurfaces follows immediately from (3.2):
$$[X,Y]=\nabla_XY-\nabla_YX=-\gamma_1\,X-\gamma_2\,Y,$$
i.e. the distribution $\Delta$ of any bi-umbilical hypersurface of type number two
is involutive.

Bi-umbilical hypersurfaces of type number two are described by the following statement.
\begin{thm}\label{T:3.1}
The integral surfaces of the distribution $\Delta$ of any bi-umbilical
hypersurface ${\mathcal M}^n$ of type number two are two-dimensional spheres.

Conversely, any two-dimensional sphere generates a family of bi-umbilical hypersurfaces
of type number two.
\end{thm}
\emph{Proof:} Let $\mathcal M^2: \; z=z(u,v), \; (u,v) \in \mathcal D$
be a fixed regular integral surface of the distribution $\Delta$. Further we assume that
the parametric lines of $\mathcal M^2$ are orthogonal and denote
$z_u^2=E, \; z_v^2=G$. Then we can choose
$\ds{X=\frac{z_u}{\sqrt E}}$, \; $\ds{Y=\frac{z_v}{\sqrt G}}$.
Taking into account (3.2), we obtain the following Frenet type formulas for
the surface $\mathcal M^2$:
$$
\begin{array}{lccl}
\nabla_XX= & \gamma_1\,Y & +  \, \ds{\sum_{i=1}^{n-2}\lambda_i \, e_i} &+\,\nu \,N,\\
[4mm]
\nabla_XY=-\gamma_1\,X, &  & &\\
[4mm]
\nabla_Xe_i=-\lambda_i\,X,\\
[4mm] \nabla_XN=-\nu\,X,\\
[6mm]
\nabla_YX= &\gamma_2\,Y,\\
[4mm]
\nabla_YY=-\gamma_2\,X & & + \,\ds{\sum_{i=1}^{n-2} \lambda_i \, e_i} &+\,\nu\,N,\\
[4mm]
\nabla_Ye_i=& -\lambda_i\,Y, & & \\
[4mm]
\nabla_YN=& -\nu \,Y. & &\\
\end{array}\leqno(3.3)
$$
In view of the equalities
$$\nabla_X\nabla_YN-\nabla_Y\nabla_XN=\nabla_{[X,Y]}N,\leqno(3.4)$$
$$\nabla_X\nabla_Ye_i-\nabla_Y\nabla_Xe_i=\nabla_{[X,Y]}e_i, \quad i=1,...,n-2\leqno(3.5)$$
we find that
$$X(\nu)=Y(\nu)=0, \quad X(\lambda_i)=Y(\lambda_i)=0.$$
Hence
$$\nu={\rm const}\neq 0, \quad \lambda_i=c_i={\rm const}.$$

If all $c_i=0$, then the space $\span\{z_u, z_v, N\}$ is a constant
three-dimensional subspace ${\R}^3$ of ${\R}^{n+1}$ and
$\mathcal M^2$ is a surface in ${\R}^3$. It follows
immediately from (3.3) that $\mathcal M^2$ lies on a sphere with radius $\ds{\frac{1}{|\nu|}}$.

If $\ds{\sum_i c_i^2=c^2>0}$, we consider the unit vector field
$b=\ds{\frac{1}{c}\,\sum_ic_i\,e_i}$. Then the equations (3.3) become
$$
\begin{array}{lccl}
\nabla_XX= & \gamma_1\,Y & + \,c \, b & + \, \nu \,N \\
[4mm]
\nabla_XY=-\gamma_1\,X, &  & &\\
[4mm]
\nabla_Xb \; = -c \, X, \\
[4mm]
\nabla_XN=-\nu \, X;\\
[6mm]
\nabla_YX= &\gamma_2\,Y,\\
[4mm]
\nabla_YY=-\gamma_2\,X & & +\,c\,b &+\,\nu\,N,\\
[4mm]
\nabla_Yb\;=& -c\,Y, & & \\
[4mm]
\nabla_YN=& -\nu \,Y. & &\\
\end{array}\leqno(3.6)$$

The equalities (3.6) imply that the space $\span\{z_u, z_v, b, N\}$ is a constant
four-dimensional subspace ${\R}^4$ of ${\R}^{n+1}$ and $\mathcal M^2$ lies in ${\R}^4$. Further
it follows that the vectors
$$\bar b=\cos \alpha \, b+\sin \alpha \, N, \quad \bar N=-\sin \alpha \, b +\cos \alpha \,N, $$
where $\tan \alpha =-\ds{\frac{c}{\nu}}$ satisfy the conditions
$$\nabla_X\bar b= \nabla_Y \bar b=0, \quad \nabla_X\bar N= -\sqrt{c^2+\nu^2}\,X,
\quad \nabla_Y\bar N=-\sqrt{c^2+\nu^2}\,Y.$$

The last equalities again show that $\mathcal M^2$ lies on a sphere with radius
$\ds{\frac{1}{\sqrt{c^2+\nu^2}}}$ in a constant three-dimensional subspace
${\R}^3$ in ${\R}^4$, orthogonal to the constant vector $b$.
\vskip 2mm
Conversely, let $S^2(r): \; z=z(u,v), \; (u,v) \in \mathcal D$ be a sphere
in ${\R}^3 \subset {\R}^{n+1}$. Taking into account the above arguments, we shall give
a construction of the bi-umbilical hypersurfaces generated by a given sphere $S^2(r)$.

Let $\bar N$ denote the normal vector field to $S^2(r)$ in ${\R}^3$ and $e$ be a constant
unit vector in ${\R}^{n-2}$ orthogonal to ${\R}^3$. Choose an orthonormal basis
$\{e, e_2,...,e_{n-2}\}$ of ${\R}^{n-2}$ and consider the vectors
$$\begin{array}{l}
N=\; \; \; \cos \alpha \, \bar N +\sin \alpha \, e,\\
[2mm]
e_1=-\sin \alpha \, \bar N + \cos \alpha \, e,
\end{array} \quad \alpha = {\rm const}.
$$
Then the required hypersurfaces are constructed as follows:
$$\mathcal M^n: \; X(u,v,w^1,...,w^{n-2})=z(u,v)+\sum_{\alpha=1}^{n-2} w^{\alpha}e_{\alpha};
\quad (u,v)\in \mathcal D, \; w^{\alpha}\in {\R}.$$
Direct calculations show that ${\mathcal M}^n$ is a hypersurface of type number two
with normal vector field $N$ and $\nu_1=\nu_2$. Hence $\mathcal M^n$
is a bi-umbilical hypersurface of type number two. \qed
\vskip 2mm
Let ${\mathcal M}^n$ be a hypersurface in $\R^{n+1}$ of type number two, which is the envelope
of a two-parameter family of hyperplanes $\{\R^n(u,v)\}, \,(u,v)\in {\mathcal D}$, defined
in a domain ${\mathcal D} \subset \R^2$ \cite{GM2}. We denote by $l = l(u, v)$ the
unit normal vector field of the hyperplane $\R^n(u,v)$ ($l$ is determined up to a sign) and by
$r = r(u, v)$ - the oriented distance from the origin of $\R^{n+1}$ to $\R^n(u,v)$. We
assume that $l$, $l_u$ and $l_v$ are linearly independent (otherwise ${\mathcal M}^n$ is
a developable ruled hypersurface, which is a hypersurface with type number one). Then,
locally there exist $n-2$ mutually orthogonal unit vectors
$b_1(u,v), \dots, b_{n-2}(u,v), \,(u,v)\in {\mathcal D}$,
which are orthogonal to $\span \{l, l_u, l_v\}$. Further we denote
$$E(u,v) = g(l_u, l_u), \quad F(u,v) = g(l_u, l_v), \quad G(u,v) = g(l_v, l_v),$$
satisfying the inequalities $E > 0, \,\, G > 0, \,\, E\,G - F^2 > 0$.

Then ${\mathcal M}^n$ can be parameterized locally as follows:
$$X(u, v, w^{\alpha}) = r\,l + \displaystyle{\frac{G\,r_u - F\,r_v}{W^2}\,\,l_u} +
\displaystyle{\frac{E\,r_v - F\,r_u}{W^2}\,\,l_v} + \sum_{\alpha =
1}^{n-2} w^{\alpha}\,b_{\alpha} , \leqno{(3.7)}$$
where $(u,v) \in {\mathcal D}, \,\, w^{\alpha} \in \R, \,\, \alpha =
1, \dots, n-2$.

Thus each pair of a unit vector-valued function $l = l(u, v)$ and a scalar function
$r = r(u, v)$ determines a hypersurface ${\mathcal M}^n$ of type number two by the equation (3.7).

In \cite{GM2} we proved that a hypersurface of type number two given by (3.7) is
bi-umbilical if and only if
$$\begin{array}{cl}
l_{uu} - l_{vv} & = \displaystyle{\frac{E_u}{E}\,l_u - \frac{E_v}{E}\,l_v},\\
[2mm]
2 l_{uv} & = \displaystyle{\frac{E_v}{E}\,l_u + \frac{E_u}{E}\,l_v},\\
[3mm]
r_{uu} - r_{vv}  & = \displaystyle{\frac{E_u}{E}\,r_u - \frac{E_v}{E}\,r_v},\\
[2mm]
2 r_{uv} & = \displaystyle{\frac{E_v}{E}\,r_u + \frac{E_u}{E}\,r_v}.
\end{array} \leqno(3.8)$$

Taking into account that the vector function $l$ is the normal vector field to
${\mathcal M}^n$, Theorem \ref{T:3.1} implies that the solutions to the system
(3.8) can be found explicitly.
\vskip 1mm
In what follows we shall consider hypersurfaces of type number two satisfying the condition
$$\nu_1(p)-\nu_2(p)\neq 0, \quad p \in {\mathcal M}^n.\leqno(3.9)$$

\subsection{Minimal hypersurfaces of type number two with involutive distribution}
Let ${\mathcal M}^n$ be a regular hypersurface of type number two satisfying the condition (3.9).
Formulas (3.1) imply that the distribution $\Delta$ is involutive if and only if
$$\sigma(x_0)=0, \quad x_0\in \Delta^{\perp}.\leqno(3.10)$$

Further we assume that ${\mathcal M}^n$ is with involutive distribution $\Delta$, i.e.
the condition (3.10) is valid.

We denote by $\mathcal K_0$ the class of hypersurfaces of type number two with involutive
distribution $\Delta$.

A hypersurface ${\mathcal M}^n$ is minimal if $\nu_1+\nu_2=0$.

Minimal hypersurfaces of type number two with involutive distribution are described
by the following statement.

\begin{thm}\label{T:4.1}
The integral surfaces of the distribution $\Delta$ of any minimal hypersurface ${\mathcal M}^n$
of the class $\mathcal K_0$ is a minimal surface in $\R^3$ or in ${\Sp}^3(r)$.

Conversely, any minimal surface in $\R^3$ or in ${\Sp}^3(r)$ generates a minimal hypersurface
of the class $\mathcal K_0$.
\end{thm}
\emph{Proof:} We put
$$\nu_1=\nu>0, \qquad \nu_2=-\nu.\leqno(3.11)$$

Let $\mathcal M^2: \; z=z(u,v), \; (u,v) \in \mathcal D$ be a fixed regular integral
surface of the distribution $\Delta$. We assume that the parametric lines of $\mathcal M^2$
are orthogonal and denote $z_u^2=E, \; z_v^2=G$. Then we consider the unit vector fields
$\ds{X=\frac{z_u}{\sqrt E}}$, \; $\ds{Y=\frac{z_v}{\sqrt G}}$.
Taking into account (3.10) and (3.11), we obtain from (3.2) the following Frenet type formulas for
the surface $\mathcal M^2$:
$$
\begin{array}{lccl}
\nabla_XX= &\; \gamma_1\,Y & + \, \ds{\sum_i \lambda_i \, e_i} &+\,\nu\,N,\\
[3mm]
\nabla_XY=-\gamma_1\,X &  & &\\
[3mm]
\nabla_Xe_i=-\lambda_i\,X & & &\\
[3mm]
\nabla_XN=-\nu\,X,\\
[6mm]
\nabla_YX= & \; \gamma_2\,Y & &\\
[3mm]
\nabla_YY=-\gamma_2\,X & & +\,\ds{\sum_i \lambda_i \, e_i} &-\,\nu \,N,\\
[3mm]
\nabla_Ye_i=&-\lambda_i \,Y, & & \\
[3mm]
\nabla_YN=& \; \nu\,Y. & &
\end{array}\leqno(3.12)$$

Using (3.12) we find from (3.5)

$$(\lambda_i)_u=0, \quad (\lambda_i)_v=0.$$

The last equalities imply that
$$\lambda_i=c_i={\rm const}, \quad i=1,...,n-2.$$

If all $c_i=0$, then the space $\span\{z_u, z_v, N\}$ is a constant
${\R}^3$ and $\mathcal M^2:\; z=z(u,v)$ is a surface in ${\R}^3$.
It follows immediately from (3.12) that $\mathcal M^2$ is a minimal surface in ${\R}^3$.

If $\ds{\sum_i c_i^2=c^2>0}$, we consider the unit vector field
$b=\ds{\frac{1}{c}\,\sum_ic_i\,e_i}$. Then (3.12) become
$$
\begin{array}{lccl}
\nabla_XX= & \gamma_1\,Y & + \, c\, b &+\,\nu \,N,\\
[4mm]
\nabla_XY=-\gamma_1\,X, &  & &\\
[4mm]
\nabla_Xb=-c \,X,\\
[4mm] \nabla_XN=-\nu\,X;\\
[6mm]
\nabla_YX= &\gamma_2\,Y,\\
[4mm]
\nabla_YY=-\gamma_2\,X & & +\,c\,b &-\,\nu\,N,\\
[4mm]
\nabla_Yb=& -c\,Y, & & \\
[4mm]
\nabla_YN=& \nu \,Y. & &\\
\end{array}$$

These equalities show that $\mathcal M^2$ is a minimal surface in
${\Sp}^3 \, (r=1/c)$.

Conversely, let $\mathcal M^2: \; z=z(u,v), \; (u,v) \in \mathcal D$ be a minimal
surface in a fixed ${\Sp}^3 (r)$ $\subset {\R}^{n+1}$. Taking into account the above
arguments, we shall give a construction of the minimal hypersurface of type number
two generated by the given minimal surface $\mathcal M^2$.

Let $N$ be the normal vector field to $\mathcal M^2$ in ${\Sp}^3$. Denote by
${\R}^{n-2}(u,v)$ the plane in ${\R}^{n+1}$ orthogonal to $\span\{z_u, z_v, N\}$ at each
point of $\mathcal M^2$ and choose a base $\{b_{\alpha}(u,v)\}, \; \alpha=1,...,n-2$ of
${\R}^{n-2}(u,v)$. Then we consider the hypersurface in ${\R}^{n+1}$ given by
$$\mathcal M^n: \; X(u,v;w^1,...,w^{n-2}):=z(u,v)+\sum w^{\alpha}b_{\alpha},
\quad (u,v) \in \mathcal D, \quad w^{\alpha}\in \R, \, \alpha=1,...,n-2.$$
By straightforward computations it follows that $\mathcal M^n$ is a minimal hypersurface
of type number two. It is clear that the distribution $\Delta$ of ${\mathcal M}^n$
is involutive.

Further, let $\mathcal M^2: \; z=z(u,v), \; (u,v) \in \mathcal D$
be a minimal surface in a fixed ${\R}^3 \subset {\R}^{n+1}$. Choose an orthonormal basis
$\{b_{\alpha}\}, \; \alpha=1,...,n-2$ for the orthogonal complement ${\R}^{n-2}$ of
${\R}^3$ in ${\R}^{n+1}$ and consider the hypersurface in ${\R}^{n+1}$ given by
$$\mathcal M^n: \; X(u,v;w^1,...,w^{n-2}):=z(u,v)+\sum w^{\alpha}b_{\alpha},
\quad (u,v) \in \mathcal D, \quad w^{\alpha}\in \R, \, \alpha=1,...,n-2.$$
It is easy to check that $\mathcal M^n$ is a minimal hypersurface of type number two
with involutive distribution $\Delta$. \qed
\vskip 2mm
Let ${\mathcal M}^n$ be a hypersurface of type number two given by (3.5). In \cite{GM2}
we proved that ${\mathcal M}^n$ is minimal if and only if
$$\begin{array}{l}
l_{uu} + l_{vv} + 2E\,l = 0,\\
[2mm]
r_{uu} + r_{vv} + 2E\,r = 0,
\end{array} \qquad l_u^2=l_v^2=E, \quad l_u \, l_v=0.$$

Theorem \ref{T:4.1} means that the solutions to the above system are generated by the minimal
surfaces in $\R^3$ and by the minimal surfaces in ${\Sp}^3$.

In a next paper we will show in details that the minimal ruled hypersurfaces in ${\R}^{n+1}$
are generated by the helicoids in ${\R}^3$ or by the minimal "ruled" surfaces \cite{L} in ${\Sp}^3$.


\begin{thebibliography}{99}
\bibitem{A1}
Aumann G. \emph{Die Minimalhyperregelfl\"achen.} Monatshefte f\"ur Math., \textbf{34} (1981),
293-304.

\bibitem{A2}
Aumann G. \emph{Zur Theorie Verallgemeinerter Torsaler Strahlfl\"achen.} Monatshefte f\"ur
Math., \textbf{91} (1981), 171-179.

\bibitem{L}
Lawson, H. B. \emph{Complete minimal surfaces in S3}, Ann. Math., \textbf{92} (1970), 335-374.

\bibitem{BKV}
Boeckx E., O. Kowalski and L. Vanhecke. \emph{Riemannian manifolds of conullity two.}
Singapore, World Scientific, 1996.

\bibitem{FG}
Frank H. and O. Giering. \emph{Verallgemeinerte Regelfl\"achen.} Math. Z. \textbf{150}
(1976), 261-271.

\bibitem {GM}
Ganchev G. and V. Mihova. \emph{On the Invariant Theory of Weingarten Surfaces in
Euclidean Space}. arXiv:0802.2191v1 [math.DG]

\bibitem{GM1}
Ganchev G. and V. Milousheva. \emph{Foliated semi-symmetric hypersurfaces in Euclidean space
with involutive geometric two-dimensional distribution.} Compt. Rend. Acad. Bulg. Sci.,
\textbf{59} (2006), 1, 5-10.

\bibitem{GM2}
Ganchev G. and V. Milousheva. \emph{An analytic characterization of the minimal
and the bi-umbilical foliated semi-symmetric hypersurfaces in euclidean space.}
Compt. Rend. Acad. Bulg. Sci., \textbf{60} (2007).

\bibitem{GM3}
Ganchev G. and V. Milousheva. \emph{Minimal Surfaces in the Four-Dimensional Euclidean Space.}
arXiv:0806.3334v1 [math.DG].

\bibitem{S}
Szab\'o Z. \emph{Structure theorems on Riemannian spaces satisfying R(X,Y).R=0, I.}
J. Differ. Geom. \textbf{17} (1982), 531-582.

\end{thebibliography}
\end{document}